\documentclass{amsart}
\usepackage{amssymb}
\usepackage{mathrsfs}
\usepackage{amsmath}
\usepackage{amsthm}

\usepackage[english]{babel}

\newtheorem{theorem}{Theorem}[section]
\newtheorem{lemma}[theorem]{Lemma}
\newtheorem{e-proposition}[theorem]{Proposition}

\newtheorem{e-definition}[theorem]{Definition\rm}
\newtheorem{remark}{\it Remark\/}
\newtheorem{example}{\it Example\/}
\newtheorem*{proof1}{Proof of Theorem \ref{thm:stand_twist}}

\newcommand{\hk}{hyperk\"{a}hler }

\newcommand{\kahl}{K\"{a}hler }
\newcommand{\ktiposp}{$K3^{[2]}$ type }

\newcommand{\kntipo}{$K3^{[n]}$ type}
\newcommand{\kntiposp}{$K3^{[n]}$ type }
\newcommand{\ie}{i.~e.~}

\begin{document}

\title{On natural deformations of symplectic automorphisms of manifolds of \kntipo}
\author{Giovanni Mongardi},
\address{Mathematisches Institut der Universit\"{a}t Bonn\\ Endenicher Allee, 60\\ 53115 Bonn, Germany}

\maketitle

\begin{abstract}
In the present paper we prove that finite symplectic groups of automorphisms of manifolds of \kntiposp can be obtained by deforming natural morphisms arising from $K3$ surfaces if and only if they satisfy a certain numerical condition.
\end{abstract}



\section{Introduction}

\label{}
The present paper is devoted to a natural question concerning deformations of automorphisms of \hk manifolds. Roughly speaking, given a $K3$ surface $S$ the group $\mathrm{Aut}(S)$ induces automorphisms of the Hilbert scheme $S^{[n]}$ of $n$ points of $S$. These automorphisms are called natural. Let $X$ be a \hk manifold deformation equivalent to some $S^{[n]}$ and let $G$ be a group of automorphisms of $X$. One can ask whether it is possible to deform $X$ together with $G$ to some $(S^{[n]},G)$, where $G$ is a group of natural automorphisms.
In the following we give a positive answer for all finite symplectic automorphism groups whose action on $H^2(X)$ is the natural one and for several different dimensions (cf. 
Theorem \ref{thm:stand_twist}). We remark that having the natural action on $H^2(X)$ is a necessary condition, since this action is constant under smooth deformations.\\
 There have been several works concerning automorphisms of $K3$ surfaces, we will refer to the foundational work of Nikulin \cite{nik1}, later improved by Mukai \cite{muk} in the nonabelian case.
By the work of Mukai \cite{muk} there are $79$ possible finite groups of symplectic automorphisms on $K3$ surfaces and, by a recent classification due to Hashimoto \cite{hashi}, there are $84$ different possibilities for their action on $H^2$. Our result holds for all these $84$ cases as long as the hypothesis of the global Torelli theorem are satisfied.\\
In the case of manifolds of \kntiposp the notion of natural morphisms was introduced by Boissi\`{e}re \cite{boi} and further analyzed by him and Sarti \cite{bs}. In the particular case of symplectic involutions on manifolds of \ktiposp our result is proven in \cite{me1}.
\subsection*{Notations}
If $L$ is a lattice and $G\subset O(L)$ we denote by $T_G(L):=L^G$ the invariant sublattice and by $S_G(L):=T_G(L)^\perp$ the coinvariant sublattice.
For $G\subset \mathrm{Aut}(X)$ and $H^2(X,\mathbb{Z})$ endowed with a quadratic form, we denote $T_G(X):=T_G(H^2(X,\mathbb{Z}))$ the invariant sublattice and $S_G(X):=S_G(H^2(X,\mathbb{Z}))$ the coinvariant sublattice. 
Let $X$ be a \hk manifold and let $G\subset \mathrm{Aut}(X)$. The group $G$ is called symplectic if it acts trivially on $H^{2,0}(X)$, \ie it preserves the symplectic form. We denote by $\mathrm{Aut}_s(X)$ the subgroup of automorphisms of $X$ preserving the symplectic form.
We will call manifolds of \kntiposp all manifold deformation equivalent to the Hilbert scheme of $n$ points on a K3 surface.

\subsection*{Preliminaries}
In this section we gather some useful results for ease of reference. The reader interested in \hk manifolds can consult \cite{huy} and \cite{huy_tor} for further references and for a broader treatment of the subject.




A \hk manifold is a simply connected compact \kahl manifold whose $H^{2,0}$ is generated by a symplectic form.
\begin{theorem}
Let $X$ be a \hk manifold of dimension $2n$. Then there exists a canonically defined pairing $(\,,\,)_X$ on $H^2(X,\mathbb{C})$, the Beauville-Bogomolov pairing, which is a deformation and birational invariant. This form makes $H^2(X,\mathbb{Z})$ a lattice of signature $(3,b_2(X)-3)$. 
\end{theorem}
For every \hk manifold $X$ and every \kahl class $\omega$ there exists a family of smooth deformations of $X$ over the base $\mathbb{P}^1$. This family is called \emph{twistor family} and denoted $TW_{\omega}(X)$.
\begin{example}
Let $X$ be a \hk manifold of \kntipo. Then $H^2(X,\mathbb{Z})$ endowed with its Beauville-Bogomolov pairing is isomorphic to the lattice 
\begin{equation}\label{latticeK3n}
L_n:=H^2(K3,\mathbb{Z})\oplus (2-2n).
\end{equation}

\end{example}
If $X$ is \hk we call a marking of $X$ any isometry between $H^2(X,\mathbb{Z})$ and a lattice $M$. There exists a moduli space of marked \hk manifolds with $H^2(X,\mathbb{Z})\cong M$ and we denote it by $\mathcal{M}_M$.\\
We will often consider the induced action of $\mathrm{Aut}(X)$ on $O(H^2(X,\mathbb{Z}))$ for a manifold $X$ of \kntipo. For a general \hk manifold this map might not be injective but in our case it is: 
\begin{lemma}\label{lem:no_identity}
Let $X$ be a manifold of \kntipo. Then the map 
\begin{equation}
\nu(X)\,:\,\mathrm{Aut}(X)\rightarrow O(H^2(X,\mathbb{Z}))
\end{equation}
 is injective.
\begin{proof}
By \cite[Theorem 2.1]{ht3} the kernel of $\nu(X)$ is invariant under smooth deformations. Beauville \cite[Lemma 3]{beau2} proved that, if $S$ is a $K3$ surface with no nontrivial automorphisms, then $\mathrm{Aut}(S^{[n]})=\{Id\}$, therefore $\{Id\}=Ker(\nu(S^{[n]}))=Ker(\nu(X))$.
\end{proof}

\end{lemma}

The following is a very important theorem which is essential in the proof of our main result. The only truly restrictive hypothesis of Theorem \ref{thm:stand_twist} is one of the hypotheses of the following:

\begin{theorem}[Global Torelli, Verbitsky, Markman and Huybrechts]\label{thm:global_torelli}
Let $X$ and $Y$ be two \hk manifolds of \kntiposp and let $n-1$ be a prime power. Suppose $\psi\,:\,H^2(X,\mathbb{Z})\,\rightarrow\,H^2(Y,\mathbb{Z})$ is an isometry preserving the Hodge structure. Then there exists a birational map $\phi\,:\,X\,\dashrightarrow\,Y$.
\end{theorem}

Let $M$ be a lattice of signature $(3,r)$. We define $\Omega_M=\mathbb{P}(\{x\in M\otimes\mathbb{C}\,\,|\,\,x^2=0,(x,\overline{x})>0\})$ as the period domain for the lattice $M$. It is an open subset of a quadric hypersurface inside $\mathbb{P}(M\otimes\mathbb{C})$.

In the particular case where $M\cong H^2(X,\mathbb{Z})$ for some \hk manifold $X$, there exists a natural map, the period map $\mathcal{P}$, between the moduli space $\mathcal{M}_M$ and the period domain $\Omega_M$.

Moreover, when Theorem \ref{thm:global_torelli} holds, two marked manifolds having the same period are birational.

The images of twistor families in $\mathcal{M}_M$ through the period map are called \emph{twistor lines}.
A fundamental property of period domains is that they are connected by twistor lines (see \cite[Proposition 3.7]{huy_tor} or \cite{beau_k3}).

\section{Deformations of pairs}
\begin{e-definition}
Let $X$ be a manifold and let $G\subset \mathrm{Aut}(X)$. A $G$ deformation of $X$ (or a deformation of the pair $(X,G)$) consists of the following data:
\begin{itemize}
\item A flat family $\mathcal{X}\rightarrow B$, $B$ connected and $\mathcal{X}$ smooth, and a distinguished point $0\in B$ such that $\mathcal{X}_0\cong X$.
\item A faithful action of the group $G$ on $\mathcal{X}$ inducing fibrewise faithful actions of $G$. 
\end{itemize}
Two pairs $(X,G)$ and $(Y,H)$ are deformation equivalent if $(Y,H)$ lies in a $G$ deformation of $X$.
\end{e-definition}
The first interesting remark is that, to some extent, all symplectic automorphism groups of a \hk manifold can be deformed:\\
\begin{remark}\label{oss:twistor_deform}
Let $X$ be a \hk manifold such that $G\subset \mathrm{Aut}_s(X)$ and $|G|<\infty$. Let $\omega$ be a $G$ invariant \kahl class. Then $TW_{\omega}(X)$ is a $G$ deformation of $X$ over $\mathbb{P}^1$.
\end{remark}
There is also a notion of local universal  $G$ deformation, for a proof of its existence we refer to \cite{me1}.
\begin{lemma}
Let $X$ be a manifold of \kntiposp and let $G\subset \mathrm{Aut}_s(X)$. Then there exists a universal local $G$ deformation of $X$ sitting inside $Def(X)$. It is locally given by the $G$-invariant part of $H^1(T_X)$ and it is of dimension $\mathrm{rank}(T_G(X))-2$. Moreover two birational manifolds with isomorphic actions of $G$ on cohomology have intersecting local $G$-deformations.
\begin{proof}
Let $X$ be birational to $Y$ and let the action of $G$ on $H^2(X)$ coincide with the action of $G$ on $H^2(Y)$ induced by the birational transformation between $X$ and $Y$. Let us take a representative $U$ of $Def(X)$ and let $x$ be a very general point inside $U^G$, which is a representative of the local $G$ deformations of $X$ and $Y$. Let $\mathcal{Y}_x$ and $\mathcal{X}_x$ be the two \hk manifolds corresponding to $x$ on $U^G$. We have $Pic(\mathcal{Y}_x)=Pic(\mathcal{X}_x)=S_G(X)$ and $\mathcal{Y}_x$ is birational to $\mathcal{X}_x$. However any $G$ invariant \kahl class on $\mathcal{Y}_x$ is orthogonal to $Pic(\mathcal{Y}_x)$ and therefore also to the set of effective curves on $\mathcal{Y}_x$, which is therefore empty. Thus the \kahl cone of $\mathcal{Y}_x$ coincides with the positive cone and $\mathcal{Y}_x=\mathcal{X}_x$. 
\end{proof}
\end{lemma}
We remark that the local $G$ deformations around two birational manifolds might not meet for a nonsymplectic group $G$.
\begin{e-definition}
Let $S$ be a $K3$ surface and let $G\subset \mathrm{Aut}_s(S)$ be a group of symplectic automorphisms on $S$. $G$ induces a group of symplectic morphisms on $S^{[n]}$ which we still denote as $G$. We call the pair $(S^{[n]},G)$ a natural pair, following \cite{boi}.
We call standard any pair $(X,H)$ deformation equivalent to a natural pair. 
\end{e-definition}
A natural question is asking under which condition a pair $(X,G)$ is standard. In the rest of the paper we make the following assumption and we prove that it is equivalent to $(X,G)$ being standard.
\begin{e-definition}
Let $X$ be a manifold of \kntiposp and let $G\subset \mathrm{Aut}_s(X)$. The group $G$ is \emph{numerically standard} if the following holds
\begin{itemize}
\item $S_G(X)\cong S_H(S)$,
\item $T_G(X)\cong T_H(S)\oplus \langle t\rangle$.
\item $t^2=-2(n-1)$, $(t,H^2(X,\mathbb{Z}))=2(n-1)\mathbb{Z}$.
\end{itemize}
For some $K3$ surface $S$ and some $H\subset \mathrm{Aut}_s(S)$ such that $H\cong G$.
\end{e-definition}
Notice that for a standard pair $(X,G)$ the group $G$ is numerically standard, since by \cite{bs} a natural pair is numerically standard.
Now the main result of the paper can be explicitly stated:
\begin{theorem}\label{thm:stand_twist}
Let $X$ be a manifold of \kntiposp and let $n-1$ be a prime power. Let $G\subset \mathrm{Aut}_s(X)$ be a finite group of numerically standard automorphisms. Then $(X,G)$ is a standard pair.
\end{theorem}

In this section we prove Theorem \ref{thm:stand_twist} using some properties of a particular period domain defined by the action of a finite group $G$ of symplectic automorphisms of a manifold $X$ of \kntipo.
\begin{e-definition}
Let $M$ be a lattice of signature $(3,r)$ and let $G\subset O(M)$. We call $\Omega_{G,M}$ the set of points $\omega$ in the period domain $\Omega_M$ such that $\omega\,\in\,T_G(M)\otimes\mathbb{C}$. 
\end{e-definition}
\begin{e-definition}
Let $\mathcal{M}_n:=\mathcal{M}_{L_n}$ be the moduli space of marked manifolds of \kntiposp and let $G\subset \mathrm{Aut}_s(X)$ for some marked $(X,f)\in\mathcal{M}_n$. Let us denote with $G$ the group of isometries induced by $G$ on the lattice $L_n$ and let $\Omega_{G,n}:=\Omega_{G,L_n}$ be as above. Then we define $\mathcal{M}_{G,n}\subset\mathcal{M}_n$ as the counterimage through the period map of $\Omega_{G,n}$.
\end{e-definition}
By the following remark the set $\mathcal{M}_{G,n}$ is the set of marked pairs $(X,f)$ such that $f^{-1}(S_G(L_n))\subset Pic(X)$ for an appropriate marking $f$ and $\Omega_{G,n}$ is just the period domain $\Omega_{T_G(L_n)}$.
\begin{remark}
Let $X$ be a \hk manifold and let $G\subset \mathrm{Aut}_s(X)$ be a finite group. Then $T_G(X)$ contains $T(X)$ and $S_G(X)\subset Pic(X)$. Moreover $T_G(X)$ has signature $(3,r)$ for some $r\geq 0$. A proof of this fact can be found in \cite[Proposition 6]{beau2}.

\end{remark}

This means that, through a chain of twistor families, we can connect any marked point $(X,f)\in\mathcal{M}_{G,n}$ with $G\subset \mathrm{Aut}_s(X)$ numerically standard to a marked point $(Y,g)$ that has the same period of a natural pair $(S^{[n]},G)$ for an appropriate marking $f'$ of $S^{[n]}$.
Since by Remark \ref{oss:twistor_deform} twistor families are $G$ deformations, we have that $(X,G)$ and $(Y,G)$ are deformation equivalent. 
\begin{proof1}
Let $X$ be a manifold of \kntiposp and let $n-1$ be a prime power. Let $G\subset \mathrm{Aut}_s(X)$ be a finite numerically standard group of symplectic automorphisms. Since $\Omega_{G,n}$ is connected by twistor lines, $(X,G)$ is deformation equivalent to $(Y,G)$ and $\mathcal{P}(Y,f)=\mathcal{P}(S^{[n]},f')\in\Omega_{G,n}$. Here $S$ is a $K3$ surface with $G\subset \mathrm{Aut}_s(S)$ and $Pic(S)=S_G(S)$, \ie the very general $K3$ surface with $G\subset Aut_s(S)$. By Theorem \ref{thm:global_torelli} there is a birational map $\phi$ between $Y$ and $S^{[n]}$ which gives an induced action of $G$ on $S^{[n]}$ (possibly nonregular). Let us denote by $H$ the group induced on $S^{[n]}$ by $\phi$ and let us keep calling $G$ the group induced by the automorphisms of $S$. We obtain our claim by proving that $H=G$ (as actions on $S^{[n]}$), since in that case $(Y,G)$ and $(S^{[n]},H)$ would be deformation equivalent through their local universal $G$-deformations.\\
Notice that, by the assumption on the numerical standardness, the actions of $G$ and $H$ already coincide on $H^2(S^{[n]},\mathbb{Z})$. 
Let now $g\in G$ and let $h$ be the element of $H$ such that $g^*=h^*$ in $H^2(S^{[n]},\mathbb{Z})$. Let $r$ be the order of $g$. Then $g\circ h^{r-1}$ induces the identity on $H^2(S^{[n]},\mathbb{Z})$. Therefore, by Lemma \ref{lem:no_identity}, $g^{-1}=h^{r-1}$, which implies $G=H$ as group of automorphisms of $S^{[n]}$.
\end{proof1}

\section*{Acknowledgements}
The present paper is an improvement of a result contained in my PhD thesis. I am very grateful to my former advisor, Prof. K. G. O'Grady for his support and to Prof. D. Huybrechts for his advice.\\
Partially supported by the Research Network Program GDRE-GRIFGA.

\end{document}